\documentstyle[12pt]{amsart}

\newcommand{\Adele}{{\Bbb{A}}}

\newcommand{\R}{{\Bbb R}}
\newcommand{\Q}{{\Bbb Q}}
\newcommand{\Z}{{\Bbb Z}}

\input amssym.def
\input amssym.tex
\def\emph#1{{\em #1}}
\def\tr{\mathrm{tr}}
\def\GL{\mathrm{GL}}
\def\ellip{{\mathrm{ell}}}
\newcommand{\mathfrak}{\frak}
\newtheorem{thm}{Theorem}
\newtheorem{Thm}{Theorem}[subsection]
\newtheorem{Lem}[Thm]{Lemma}

\newtheorem{Prop}[Thm]{Proposition}
 
\theoremstyle{definition}

\newtheorem{Rem}{Remark}   
\newtheorem{Fact}{Fact}
\numberwithin{equation}{subsection}  

\renewcommand{\rm}{\normalshape}
\newcounter{SecNumb}

\title[Lecture \FOO]
{
Lectures on the Arthur--Selberg Trace Formula
}
\author{Steve Gelbart}
\date{MSRI, Spring 1995}
\begin{document}
\maketitle

These are Notes prepared for nine lectures given at the
Mathematical Sciences Research Institute, 
MSRI, Berkeley during the period January--March 1995. It is a pleasant duty
to record here my gratitude to MSRI, and its staff, for making possible this
1994--95 Special Year in Automorphic Forms, and for providing such a setting
for work.

The purpose of these Notes is to describe the contents of Arthur's
earlier,  foundational papers on the trace formula.  In keeping with 
the introductory nature
of the lectures, we have sometimes illustrated the ideas of 
Arthur's general theory by
applying them in detail to the case of $\GL(2)$;  we have also included
a few lectures on the ``simple trace formula'' (and its applications), 
and on Jacquet's relative trace formula.

The {\TeX} preparation of these Notes I owe to Wendy McKay, who patiently
and professionally transformed my messy scrawl into something readable; 
her expertise, and good cheer under the pressure of weekly deadlines, is
something I shall not soon forget.

I  wish to thank the auditors of these lectures for their interest, and ,
  J.~Bernstein, D.~Goldberg, E.~Lapid, C.~Rader, S.~Rallis,
A.~Reznikov, and D.~Soudry, for their helpful suggestions and explanations. 
Finally, I wish to thank J.~Arthur, H.~Jacquet, and J.~Rogawski for
many tutorials on these and related topics over the past year;  I hope they do
not mind seeing some of their comments reappear in these Notes.

\vskip 18pt
\section*{Table of Contents}

{
\obeylines\parskip=.5\baselineskip
Lecture I.  Introduction to the Trace Formula

Lecture II.  Arthur's Modified Kernels I:  The Geometric Terms

Lecture III.  Arthur's Modified Kernels II: The Spectral Terms

Lecture IV.  More Explicit Forms of the Trace Formula

Lecture V.  Simple Forms of the Trace Formula

Lecture VI.  Applications of the Trace Formula

Lecture VII.  $(G,M)$-Families and the Spectral $J_\chi(f)$

Lecture VIII.  Jacquet's Relative Trace Formula

Lecture IX.  Some Applications of Paley--Wiener, and Concluding Remarks

}
\vfil\eject
\addtocounter{section}1

\null
\vspace{30pt}
\section*{Lecture I. Introduction to the Trace Formula}
\vspace{30pt}

\def\FOO{I. INTRODUCTION TO THE TRACE FORMULA}

The ultimate purpose of the trace formula is to better understand the structure
of automorphic forms on a general reductive group $G$. In particular, one wants
to understand the decomposition of the right regular representation $R_0(g)$
of $G(\Adele)$ in the space of cusp forms $L^2_0(G(F)\setminus G(\Adele))$.
The original idea of the (Selberg) trace formula---roughly speaking---was to
describe the trace of the integral operator $R_0(f)=\int_{G(\Adele)}f(g)R_0(g)dg$ for $f$ in $C_c^\infty (G(\Adele))$; this Selberg
succeeded to do for certain special $G$.
On the other hand, \emph{Arthur's trace formula} is valid for arbitrary $G$, and asserts the equality of two
different sums of distributions on $G$, namely:
\begin{equation}
 \sum_{\mathfrak{o}\in\mathfrak{O}}J_{\mathfrak{o}}^T(f)=\sum_{\chi\in\mathfrak{X}}J_\chi^T(f).
\label{eq:Arthur}
\end{equation}
Here the left side of (\ref{eq:Arthur}) is summed over certain equivalence classes in the group
of rational points $G(F)$; it is called the \emph{geometric side} of the trace
formula. On the right side of (\ref{eq:Arthur}), the sum is over certain ``cuspidal data''
$\chi=\{(M,\sigma)\}$ of $G$; this is the \emph{spectral side} of the trace formula (one of whose terms is $\tr R_0(f)$).

Our purpose in these lectures is to describe and develop this formula in complete detail, but mostly in the context of $\GL(2)$. We also wish to describe some
refinements of (\ref{eq:Arthur}) necessary for applications to automorphic forms.
The basic references are Arthur's original papers [A1]--[A10]. Earlier
expositions of the $\GL(2)$ trace formula, such as [Ge] and [GJ], will be referred to occasionally; but as they mostly predate Arthur's general development of the trace formula, their \emph{point of view} will largely be ignored, and their
arguments not reproduced here.

\bigskip
\subsection{Some history}

``The'' trace formula was introduced by Selberg, in the context of
a semi-simple Lie group $G$ and discrete subgroup $\Gamma$, in his famous
1956 paper [Se]. In that paper Selberg first of all described a general formula in
the case of \emph{compact quotient} $\Gamma\setminus G$ (see the next section below
for a brief description). Secondly, he treated in detail certain special
\emph{non-compact} quotient cases such as $SL_2(\Z)\setminus SL_2(\R)$; here
Selberg analytically continued Eisenstein Series $E(z,s)$ in order to handle
the continuous spectrum of $L^2(\Gamma\setminus G)$, and derived a formula for the
trace of $R(f)$ in the space of cusp forms.
In both the case of compact and non-compact quotient, Selberg was motivated
by applications to geometry and number theory (lengths of geodesics in the
Riemann surface in the case of compact quotient, traces of Hecke operators
in case of non-compact quotient, etc.). But despite the obvious
desirability of generalizing these results in the case of non-compact
quotient, the obstacles---such as analytically continuing the Eisenstein
series in higher rank---were formidable.
A big breakthrough came in the 1960's when Langlands developed a theory of
Eisenstein series valid for any reductive group $G$, and used it to describe the continuous spectrum of $L^2(G(F)\setminus G(\Adele))$ (cf. [La 1]).
As an outgrowth of this theory, Langlands was led to his amazing conjectures
about automorphic forms and their $L$-functions, i.e., the ``Langlands
program'', with its cornerstone the \emph{principle of functoriality of automorphic forms}. A typical result in this program is the following:

\begin{thm}[{[JL]}]
Let $D$ be a division quaternion algebra over a number field $F$ and $G'=D^*$. Then to each automorphic cuspidal representation $\sigma$
of $G'$ (i.e., an irreducible $G'(\Adele)$-submodule of $L^2(Z(\Adele)G'(F)\setminus G'(\Adele))$ of dimension greater than one) there exists a corresponding automorphic
cuspidal representation $\pi=\pi(\sigma)$ of $G(\Adele)=\GL_2(\Adele)$ (an
irreducible $G(\Adele)$-submodule of the space of cusp forms $L^2_0(Z(\Adele)G(F)\setminus G(\Adele))$ with the property that for each place $v$ in $F$ unramified for $D$
(i.e. where $G'(F_v)\simeq \GL_2(F_v)$), $\sigma_v\simeq \pi_v$.
\end{thm}

Interestingly enough, this correspondence (the well-known
\emph{Jacquet--Langlands correspondence}) can be proved using any one of the three basic tools of
the theory of automorphic forms: the theory of automorphic \emph{$L$-functions},
the theory of \emph{theta-series} liftings, \emph{or} the trace formula (for
both $G'$ and $G=\GL_2$, essentially as developed by Selberg).
However, in this case at least, the trace formula approach gives the strongest
payoff, namely a \emph{characterization} of the image of the correspondence
$\sigma\rightarrow\pi(\sigma)$: a cuspidal representation $\pi$ of $G(\Adele)$
is of the form $\pi(\sigma)$ for some $\sigma$ on a given $G'_D$ if and only if
for each \emph{ramified} prime $v$, $\pi_v$ is a square integrable (discrete
series) representation of $G_v=\GL_2(F_v)$.

Subsequently, a (``twisted'') trace formula for $\GL_2$ 
was developed by Saito-Shintani and Langlands to prove another important 
instance of functoriality,
namely \emph{base change} 
(with spectacular applications to Artin's conjecture
for Galois representations, 
and hence ultimately to Wiles' proof of Fermat\dots).
Thus the need to generalize the 
trace formula to a general reductive group $G$
was clear, and this is precisely 
what Arthur did in his papers between the
early 1970's and late 1980's.

\bigskip
\subsection{The Case of Compact Quotient}

In order to motivate Arthur's work, let us explain how the form (2) takes in
the case of compact quotient. Thus $G$ is a semi-simple group over (the number
field) $F$, with $G(F)\setminus G(\Adele)$ compact.
In particular, $G$ is anisotropic, without unipotent elements.
With $R(g)$ the right regular representation in $L^2(G(F)\setminus G(\Adele))$ and $f$
in $C_c^\infty(G(\Adele))$, a straightforward computation gives
\[(R(f)\phi)(x)=\int_{G(\Adele)}f(y)\phi(xy)\,dy=\]
\[ \int_{G(\Adele)}f(x^{-1}y)\phi(y)\,dy=\int_{G(F)\setminus G(\Adele)}(\sum_{\gamma\in G(F)}f(x^{-1}\gamma y))\phi(y)\,dy,\]
i.e., $R(f)$ is an integral operator on $\phi$ in $L^2(G(F)\setminus G(\Adele))$, with
kernel
\[ K_f(x,y)=\sum_{\gamma\in G(F)}f(x^{-1}\gamma y).\]
Since $G(F)\setminus G(\Adele)$ is compact, it is then a well known theorem (cf. [GGPS],
Chapter I), that $R(f)$ is of trace class, and that
\[ \tr R(f)=\int_{G(F)\setminus G(\Adele)} K_f(x,x)\,dx,\] the integral of the kernel over the ``diagonal subgroup''
$G(F)\setminus G(\Adele)^\Delta$. In this case, \emph{the (Selberg) trace formula amounts
to a calculation of $trace R(f)$ using two different expressions for this kernel.}
On the one hand, we can write
\[ K_f(x,y)=\sum_{\mathfrak{o}\in\mathfrak{O}}K_{\mathfrak{o}}(x,y)\]
where $\mathfrak{O}$ denotes the set of conjugacy classes in $G(F)$, and
\[ K_{\mathfrak{o}}(x,y)=\sum_{\gamma\in\mathfrak{o}=\{\gamma_0\}}f(x^{-1}\gamma y)=\sum_{\delta\in G_{\gamma_0}(F)\setminus G(F)}f(x^{-1}\delta^{-1}\gamma_0\delta x)\]
where $G_{\gamma_0}(F)$ is the centralizer of $\gamma_0$ in $G(F)$. Then
\[ J_{\mathfrak{o}}(f)=\int_{G(F)\setminus G(\Adele)} K_f(x,x)\,dx=\int_{G(F)\setminus G(\Adele)} \sum_{\delta\in G_{\gamma_0}(F)\setminus G(F)}f(x^{-1}\delta^{-1}\gamma_0\delta x)\,dx=\]

\[ \int_{G_{\gamma_0}(F)\setminus G(\Adele)}f(x^{-1}\gamma_0 x)\,dx=meas(G_{\gamma_0}(F)\setminus G_{\gamma_0}(\Adele))\int_{G_{\gamma_0}(\Adele)\setminus G(\Adele)}f((x^{-1}\gamma_0 x)\,dx,\] and thus
\[ \tr R(f)=\sum_{\mathfrak{o}\in\mathfrak{O}}J_{\mathfrak{o}}(f)=\sum_{\{\gamma\}}meas(G_{\gamma}(F)\setminus G_{\gamma}(\Adele))\int_{G_\gamma(\Adele)\setminus G(\Adele)}f((x^{-1}\gamma x)\,dx.\]
On the other hand, the compactness of $G(F)\setminus G(\Adele)$ also implies a discrete
decomposition $R(g)=\sum_{\pi}m_\pi\pi(g)$ where $\pi$ runs through the
irreducible ``automorphic'' representations of $G(\Adele)$, and $m_\pi<\infty$ is
the multiplicity of $\pi$ in $R$. Equivalently, $L^2(G(F)\setminus G(\Adele))=\sum L_\pi$
where $L_\pi$ is the invariant subspace realizing the $m_\pi$ copies of $\pi$,
and $R(f)=\sum R_\pi(f)$ where each $R_\pi(f)$ is an integral operator in $L_\pi$ with kernel
\[ K_\pi(x,y)=\sum_{\phi\mathrm{\ an\ ON\ basis\ for }L_\pi}(R(f)\phi)(x)\overline{\phi(y)}.\]
Thus we also have
\[\tr R(f)=\sum_{\chi\in\{(G,\pi)\}}J_\chi(f),\] where \[J_\chi(f)=\int_{G(F)\setminus G(\Adele)} K_\pi(x,x)\,dx.\]
The power of the trace formula derives from the fact that this equality
\begin{equation}
\sum J_{\mathfrak{o}}(f)=\sum J_\chi(f),
\label{eq:compactf}
\end{equation}
namely
\[ \sum_{\{\gamma\}}m_\gamma\int_{G_\gamma(\Adele)\setminus G(\Adele)}f(x^{-1}\gamma x)\,dx=\sum_{\pi}\int_{G(F)\setminus G(\Adele)}(\sum_{\phi} (R(f)\phi)(x)\overline{\phi(x)})\,dx\]
relates two completely different kinds of data, geometric and spectral (the
latter of which we rarely know much about, but very much want to!).

\bigskip
\subsection{Difficulties in the Non-Compact Case}

Arthur's great contribution was to figure out how to overcome all the obstacles
that present themselves in generalizing (\ref{eq:compactf}) from the compact to the non-compact quotient setting.
Some of the principal problems in the non-compact setting are these:

(a) Although $R(f)$ is still an integral operator, with kernel
\[ K_f(x,y)=\sum_\gamma f(x^{-1}\gamma y),\] this kernel is \emph{no longer} integrable over the
diagonal; and

(b) The regular representation $R(g)$ no longer decomposes discretely, Eisenstein series are required to describe the continuous spectrum, and $R(f)$ is of
course no longer of trace class.
Nevertheless, Arthur derives a formula (\ref{eq:Arthur})
which still relates ``geometric'' and ``spectral'' distributions attached to $G$---and has powerful applications to number theory (even if neither side of (\ref{eq:Arthur}) represents a formula for $\tr R(f)$!).
For example, in this formula, 
\[J_{\mathfrak{o}}^T(f)=\int_{Z(\Adele)G(F)\setminus G(\Adele)}K_{\mathfrak{o}}^T(x,x)\,dx,\] where $K_{\mathfrak{o}}^T$ is a \emph{modification} of the kernel
\[ K_{\mathfrak{o}}(x,x)=\sum_{\gamma\in\mathfrak{o}}f(x^{-1}\gamma x)\] which \emph{can} be integrated over the
diagonal. To see why this modification is really needed, and what it accomplishes, let us look at the example of $G=\GL(2)$, and $\mathfrak{o}$ the hyperbolic conjugacy class
\[ \mathfrak{o}=\{\delta^{-1}\bigl(\textstyle{\alpha\atop0}{0\atop1}\bigr)\delta:\delta\in M(F)\setminus G(F)\},\]
where $M=\{\left(\begin{array}{cc}a&0\\0&b\end{array}\right)\}$. In this case,
\[ J_{\mathfrak{o}}(f)=\int_{Z(\Adele)G(F)\setminus G(\Adele)}K_{\mathfrak{o}}(x,x)\,dx\]
(where $Z$ is the center of scalar matrices).
A formal computation then shows
\begin{eqnarray*}
J_{\mathfrak{o}}(f)&=&\int_{Z(\Adele)G(F)\setminus G(\Adele)}\sum_{\delta\in M(F)\setminus G(F)}f(x^{-1}\delta^{-1}\bigl(\textstyle{\alpha\atop0}{0\atop1}\bigr)\delta x)\,dx\cr&=&
\int_{Z(\Adele)G(F)\setminus G(\Adele)}f(x^{-1}\bigl(\textstyle{\alpha\atop0}{0\atop1}\bigr) x)\,dx\cr&=&
 \int_{M(\Adele)\setminus
G(\Adele)}(\int_{Z(\Adele)M(F)\setminus M(\Adele)}\,dy)f(x^{-1}\bigl(\textstyle{\alpha\atop0}{0\atop1}\bigr) x) \,dx,
\end{eqnarray*}
i.e.,
\[ J_{\mathfrak{o}}(f)=m(Z(\Adele)M(F)\setminus M(\Adele))F_f(\bigl(\textstyle{\alpha\atop0}{0\atop1}\bigl)),
\]
with $F_f(\bigl({\alpha\atop0}{0\atop1}\bigr))$ the nice ``orbital integral'' $\int_{M(\Adele)\setminus G(\Adele)} f(x^{-1}\bigl({\alpha\atop0}{0\atop1}\bigr) x) \,dx$ of a compactly supported function $f$ in $C_c^\infty(Z(\Adele)\setminus G(\Adele))$. Thus $J_{\mathfrak{o}}(f)$ is \emph{infinite} for the simple reason that
\[ m(Z(\Adele)M(F)\setminus M(\Adele))=m(F^*\setminus\Adele^*)=\infty.\]
To compensate for this, Arthur ``modifies'' $K_{\mathfrak{o}}$ to $K_{\mathfrak{o}}^T$, using a \emph{truncation parameter} $T=(T_1,T_2)$ in the Lie algebra of $M$. As we shall see,
this modification operator (in this case) ``truncates'' $F^*\setminus\Adele^*$ to a
finite interval of length $T_2-T_1$, but turns $F_f$ into a ``weighted''
orbital integral!
In general, for an arbitrary $G$, 
it will turn out that $K_{\mathfrak{o}}$ will be
integrable ``over the diagonal'' only when $\mathfrak{o}$ intersects \emph{no} proper
parabolic subgroup $P$ of $G$. 
Of course for anisotropic $G$ this always
happens, but in general one really needs to modify $K_{\mathfrak{o}}$ (as in the above example),
prove that the modified kernels $K_{\mathfrak{o}}^T$ are indeed integrable, and then evaluate
them (as generalized weighted orbital integrals).

What about the spectral side? 
Here it turns out that $L^2(Z(\Adele)G(F)\setminus
G(\Adele))=\sum_{\chi\in\mathfrak{X}}L_\chi^2$, where each $L_\chi^2$
is a $G(\Adele)$-invariant (not necessary discretely decomposable)
submodule indexed by certain cuspidal data $\mathfrak{X}=\{(M,\sigma)\}$,
where $M$ is a Levi subgroup of a parabolic $P\subset G$, and $\sigma$ is
a cuspidal automorphic representation of $M(\Adele)$. Thus we have
$K(x,y)=\sum K_\chi(x,y)$, the spectral decomposition of $K_f(x,y)$
corresponding to (\ref{eq:Arthur}). Here it also turns out that
$K_\chi(x,y)$ will \emph{not} be integrable over the diagonal, unless the
cuspidal data $(M,\sigma)=(G,\pi)$ (i.e. $P=G$). Thus we again need to
deal  with \emph{modified} $K_\chi^T$ (which are integrable\dots),
eventually arriving at the trace formula (\ref{eq:Arthur}) by integrating
the equality
\begin{equation}
\sum_{\mathfrak{o}} K_{\mathfrak{o}}^T(x,y)=\sum_\chi K_\chi^T(x,y)
\label{eq:trace}
\end{equation}
over the diagonal.

\bigskip\noindent
{\bf Remarks.} (1) The expressions on either side of (\ref{eq:trace}) no longer (generally)
represent the kernel $K_f(x,y)=\sum f(x^{-1}\gamma y)$, and hence (\ref{eq:Arthur})
is no longer a real trace formula at all. Nevertheless, what is true is the
following.
Let $\mathfrak{X}(G)$ be the subset of cuspidal data $\chi=(M,\sigma)\in\mathfrak{X}$ with $M=G$. Then $\sum_{\chi\in\mathfrak{X}(G)}K_\chi(x,y)$ is the kernel of $R(f)$ restricted
to $L_0^2(Z(\Adele)G(F)\setminus G(\Adele))$, it is integrable over the diagonal, and (\ref{eq:Arthur})
asserts that
\begin{equation}
\tr R_0(f)=\sum_{\chi\in\mathfrak{X}(G)}J_\chi(f)=\sum_{\mathfrak{o}} J_{\mathfrak{o}}^T(f)-\sum_{\mathfrak{X}\setminus\mathfrak{X}(G)}J_\chi^T(f),
\label{eq:trace2}
\end{equation}
i.e., Arthur's trace formula may be viewed as a trace formula after all.

\bigskip\noindent
(2) In most applications of the trace formula to automorphic forms we want to \emph{compare} the formula for one group with the trace formula for another.
In that case the form of (\ref{eq:Arthur}) is no less useful than (\ref{eq:trace2}).

\bigskip
\subsection{Possible Applications}

We have already mentioned the applications to quaternion algebras and base change, both of which we plan to discuss (at the level of $\GL(2)$, not $\GL(n)$ -
which is the subject matter of [AC], using the more difficult $\GL(n)$ trace
formulas). We also plan to discuss (i) Jacquet's \emph{relative trace formula},
where the idea is to integrate the (different expressions for the) kernel
over interesting \emph{non}-diagonal subgroups of $G\times G$; and (ii)
applications of the \emph{simple trace formula} such as to the embedding of a (local) discrete series representation
of $\GL_2(F_v)$ into a (global) cuspidal representation of $\GL_2(\Adele)$.
In all these applications, one needs to further refine Arthur's first form of
the trace formula. This is because (as already hinted) the ``truncation'' i.e.,
modification process used to make $K_{\mathfrak{o}}$ or $K_\chi$ integrable introduces
various defects into the distributions $J_{\mathfrak{o}}^T(f)$ and $J_\chi^T(f)$.
For example, they may no longer be \emph{invariant} distributions (invariant
under conjugation by $G(\Adele)$), or \emph{factorizable} as the product of
local distributions over $G(F_v)$. Thus one needs to replace and rearrange the
sums of distributions appearing in (\ref{eq:Arthur}) by still more complicated but
ultimately more practical expressions.

\newpage

\addtocounter{section}1
\setcounter{subsection}0

\renewcommand{\thefootnote}{\fnsymbol{footnote}}

\newcommand{\SS}[1]{\vspace{.1in}\begin{flushleft} \S #1 \end{flushleft}}
\newcommand{\place}{{\nu}}
\newcommand{\Nfield}{{F}}
\newcommand{\e}{{\text{e}}}
\newcommand{\AutoF}{automorphic form}
\newcommand{\TF}{trace formula}
\newcommand{\irr}{irreducible}
\newcommand{\rep}{representation}
\newcommand{\gp}{group}
\newcommand{\LA}{Lie algebra}
\newcommand{\red}{reductive}
\newcommand{\op}{operator}
\newcommand{\ONB}{orthonormal basis}
\newcommand{\Hom}{\operatorname{Hom}}
\newcommand{\X}{\cal X}
\newcommand{\ie}{{i.}{e.}\ }
\newcommand{\eg}{{e.}{g.}\ }
\newcommand{\etc}{etc.\ }
\newcommand{\cf}{\normalshape\mediumseries see\ }
\newcommand{\NB}{{N.}{B.}\ }
\newcommand{\word}[1]{\quad  \text{#1} \quad}

\newcommand{\bs}{\backslash}
\newcommand{\Apts}[1]{{#1(\Adele)}}
\newcommand{\AmodF}[1]{{#1(\Nfield) \bs \Apts{#1}}}
\newcommand{\AmodZF}[3]{{\Apts{#3} {#2}(\Nfield) \bs \Apts{#1}}}
\newcommand{\FmodF}[2]{{#2(\Nfield) \bs #1(\Nfield)}}
\newcommand{\AmodA}[2]{{\Apts{#2} \bs \Apts{#1}}}
\newcommand{\EisenSp}[2]{{\Apts{Z}\Apts{#2}{#1}(\Nfield) \bs \Apts{G}}}
\newcommand{\Siegel}[1]{{{\frak S}_{#1}}}
\newcommand{\SiegelPiece}[1]{{{\cal B}_{#1}}}
\newcommand{\SiegelPieceT}[1]{{{\widetilde{\cal B}}_{#1}}}
\renewcommand{\d}{\, {d}}
\newcommand{\Cinfc}{C^\infty_c}
\newcommand{\orbit}{{\frak o}}
\newcommand{\Orbit}{{\frak O}}
\newcommand{\GLtwo}{{GL(2)}}
\newcommand{\abs}[1]{{\left| {#1} \right|}}
\newcommand{\norm}[1]{{\left\| \, {#1} \, \right\|}}
\newcommand{\set}[2]{{\left\{ \, {#1} \, \mid \, #2 \, \right\}}}

\newcommand{\GeoKer}{{k^T_\orbit}}
\def\bb#1#2#3{ \bigl({\textstyle{#1\atop 0}{#2\atop#3}}\bigr)}
\newcommand{\trunc}[1]{{\widehat \tau_{#1}}}

\newcommand{\LAa}[1]{\frak a_{#1} }
\newcommand{\WeylCh}[1]{{{\LAa{#1}}^+}}

\null
\vspace{30pt}
\stepcounter{SecNumb}
\section* {Lecture II.
{Arthur}'s Modified Kernels I: The Geometric Terms}
\vspace{30pt}

\def\FOO{II. {ARTHUR}'S MODIFIED KERNELS I: THE GEOMETRIC TERMS}

Our purpose here is to develop the geometric side of {Arthur}'s first form of the \TF
$$			
      \sum_{\orbit \in \Orbit} J^T_\orbit (f) = \sum_{\chi \in \frak X} J^T_\chi (f),
$$
for the \gp\ $G = \GLtwo$. In particular, we must introduce the modified kernels
$\GeoKer(x, f)$, and prove they are integrable.

\NB The treatment of the {\TF} in \cite{sG75} and \cite{GJ77} largely ignored the
integrability of these modified kernels.

\bigskip
\subsection{The Definition of  \protect$\bold{K_\orbit}\protect$} \label{S:sec1}

In the case of compact quotient, the geometric side of the \TF\ is summed over
ordinary $G(\Nfield)$-conjugacy classes $\orbit$ in  $G(\Nfield)$ (or in
$\FmodF{G}{Z}$, where $Z$ is the center of $G$). For
arbitrary $G$, it turns out to be best to make the following definition:
$\gamma_1$ and $\gamma_2$ in $\FmodF{G}{Z}$ are called
{\em equivalent\/}, or ``{\em conjugate}'', if their semi-simple components are
$G(\Nfield)$ conjugate in the usual sense. That is, $\gamma_1 \sim \gamma_2$
means, with the Jordan decompositions $\gamma_i = \gamma_i^s \gamma_i^u$
($i = 1,2$), that $\gamma_1^s$ is $G(\Nfield)$-conjugate to $\gamma_2^s$.
For each such equivalence class $\orbit$, and for
$f \in \Cinfc(\AmodA{G}{Z})$, set
$$			
      K_\orbit(x, y) = \sum_{\gamma \in \orbit} f(x^{-1} \gamma y),
$$
and let $\Orbit$ denote the collection of all such classes $\orbit$ in
$\AmodA{G}{Z}$.

Note that if $\gamma$ is elliptic, \ie if $\gamma$ is not conjugate in 
$G(\Nfield)$ to an element of any proper parabolic subgroup  $P(\Nfield)$,
then it is automatically semi-simple. (Its eigenvalues do not lie in  $\Nfield$,
and hence are distinct.) Thus our new notion of conjugacy reduces to
ordinary conjugacy.

On the other hand, suppose $\gamma$ is $G(\Nfield)$-conjugate to an element
$p$ in a proper parabolic subgroup  $P(\Nfield)$. For $\GLtwo$, this means
$P$ is the Borel subgroup  $B = \left\{ \bb a  x  b \right\} $, so modulo $Z(\Nfield)$,
$$
	p = p^s p^u  \word{with} p^s = \bb \alpha 0 1, \word{and} p^u = \bb 1 y 1.
$$

\vspace{.1in}
{\bf Case (a):\ $\alpha = 1$.} In this case, the class $\orbit$ containing  $\gamma$
consists of two ordinary unipotent conjugacy classes in $\FmodF{G}{Z}$:
the trivial class $\{1\}$, and the ordinary conjugacy class of any non-trivial unipotent
element $\bb 1 {x_0} 1 $ (since $ \bb \alpha 0 1 $ already conjugates
$\bb 1 {x_0} 1 $ to $\bb 1 {\alpha x_0} 1 $ ).

\vspace{.1in}
{\bf Case (b):\ $\alpha \ne 1$.} In this case, $\orbit$ is an ordinary hyperbolic
conjugacy class of the form
$$
	\orbit_\alpha = \left\{ \delta^{-1} \bb \alpha 0 1 \delta
					\mid \delta \in M(\Nfield) \bs G(\Nfield) \right\}.
$$
Note that $ M(\Nfield) $ is just the centralizer $G_{\bb \alpha 0 1}(\Nfield)$,
and that $ \bb \alpha x 1 $ also belongs to $\orbit_\alpha$ (since
$ {\bb \alpha x 1} = {\bb 1 y 1}^{-1} {\bb \alpha 0 1} {\bb 1 y 1}^{-1} $ 
where $y =  {x}/{(\alpha - 1)}$\,).

\bigskip
\subsection{The Integrability of Elliptic $K_\orbit$}
\label{S:sec2}

As suggested in Lecture I, it is a general phenomenon that for an
{\em elliptic\/}
class $\orbit$, the kernel $K_\orbit(x, y)$ is integrable over the diagonal
subset  $\AmodZF{G}{G}{Z}$. We now prove this for $G = \GLtwo$, following the
general arguments of \cite{jA74}.  The crucial step is a lemma which first
requires some preliminary notions to be recalled.

According to the Iwasawa decomposition for $G$, we have
\begin{equation}  
      G(\Adele) = N(\Adele) \, M(\Adele)\, K, \label{E:Iwasawa}
\end{equation}
where $N$ is the standard unipotent subgroup $\{ \bb 1 y 1 \}$, and $K$ is
the standard maximal compact subgroup $\prod K_\nu$ of $G(\Adele)$. If we
denote by $\LAa M$ the \LA\ of $M = \{ \bb a 0 b \}$, \ie,
$$			
      \LAa M = \Hom(X(M), \R) \approx \R^2,
$$
then we can define the familiar function
$$			
      H_M : M(\Adele) \to \LAa M \word{by}
      		H_M {\bb a 0 b} = (\log \abs{a}_\Adele, \log \abs{b}_\Adele),
$$
and extend it to $G(\Adele)$ by way of the decomposition \ref{E:Iwasawa}, \ie
if $g = nak$ in $G(\Adele)$, then
$$			
      H_M(g) = H_M(a).
$$

\begin{Lem} [ {\rm \cf {\cite[p~362]{jA74}}} ] 	 \label{L:compacity}
	Suppose $\Omega$ is a subset of $G(\Adele)$ which is compact modulo
	$Z(\Adele)$. Then there exists a number $d_\Omega > 0$ with the
	following property: if $\gamma \in G(\Nfield)$ is such that
$$			
      		g^{-1} \gamma n g \in \Omega
$$
	for some $n \in N(\Adele)$ and $g \in G(\Adele)$ with
	$\alpha H(g) > d_\Omega $, then $\gamma \in B(\Nfield)$.
\end{Lem}

Assuming this Lemma (for the moment), let us quickly dispose of the
integrability of ${K_\orbit}$. For this, and for other purposes later on,
we also need to recall the notion of a {\em Siegel domain\/}, together with an
explicit integration formula for its volume.

First we refine the Iwasawa decomposition \ref{E:Iwasawa} by writing,
for each $g$ in $G(\Adele)$,
\begin{equation} \label{E:Iwasawa2}
      g = z n h_t m k
\end{equation}
where $z \in Z(\Adele)$,  $n \in N(\Adele)$, $k \in K$, 
$$			
      	h_t = {\bb {\e^t}{0}{\e^{-t}}},
$$
with $\e^t$ an idele equal to $1$ at every finite place and the real number
$\e^t$ at each infinite place, and
$$			
      	m = {\bb a 0 1},
$$
with $a \in \Adele^1 = \{a \in \Adele^\times \mid \abs{a}_\Adele = 1\}$.
The corresponding integration formula, for $f$ on $\AmodA{G}{Z}$, is
\begin{equation} 
      \int_\AmodA{G}{Z} f(g) \d g = \int\!\int\!\int\!\int_{-\infty}^{\infty}
      						f(n a t k) {\e^{-2t}} \d t \d^\times a \d n \d k. \label{E:IwasawaInt}
\end{equation}

Now for each compact subset $C_1$ of $N(\Adele)$, compact subset
$C_2$ of $\Adele^1$, and real number  $c$ (not necessarily positive),
let $\Siegel c$ denote the set of points $g$ in $G(\Adele)$ of the
form (\ref{E:Iwasawa2}) with $n \in C_1$, $a \in C_2$ and
$t > c/2$. Such a set is called a {\em Siegel set\/}. Note that
$$			
     	\Siegel {c} \subset \set{g \in \Apts{G}}{H(g) > c}.
$$
\begin{Fact}  [{\cf \cite[Thm~9]{rG63}}]		 \label{F:FundDom1}
  For any $c_1$, $c_2$,
$$			
      	\Siegel {c_1} \cap \gamma \Siegel {c_2} \ne \emptyset
$$
\noindent for only finitely many $\gamma$ in $G(\Nfield)$ {\em modulo\/} $B(\Nfield)$.
\end{Fact}
\begin{Fact}  [{\cf \cite[p~16]{rG63}}]		  \label{F:FundDom2}
If $\Siegel {c}$ is sufficiently large
(\ie $C_1$ and $C_2$ are sufficiently large, and $c$ is sufficiently
{\em small\/}), then
\begin{equation} 
      	\Apts{G} = G(\Nfield) \cdot \Siegel {c}. \label{E:GfSiegel}
\end{equation}
\end{Fact}
\setcounter{Fact}{0}
In other words, $\Siegel {c}$ is essentially a fundamental domain for
$G(\Nfield)$ acting on $\Apts{G}$; more precisely, Fact~\ref{F:FundDom2}
says only that certain $\Siegel {c}$ {\em contain\/} a fundamental domain.
Henceforth, we always assume that $\Siegel {c}$ is so chosen that
\ref{E:GfSiegel} holds, and we refer to such a $\Siegel {c}$ as a
{\em Siegel domain\/}.

\begin{Rem} From the above, it is clear that
$$			
      	\AmodZF{G}{G}{Z}
$$
has {\em finite\/} volume. In fact, by \ref{E:IwasawaInt} and \ref{E:GfSiegel},
$$			
      	\int_{\AmodZF{G}{G}{Z}} \d g \le \int_{\Apts{Z} \bs \Siegel {c}} \d g =
      	\int\!\int\!\int \left( \int_c^\infty \e^{-2t} \d t \right)
      	\d^\times a \d n \d k,
$$
with the outer three integrals over compact domains.
\end{Rem}

\begin{Prop}  \label{P:elliptic}
If $\orbit$ is {\em elliptic\/}, then
$$			
      	x \mapsto K_\orbit(x, x) = \sum_{\gamma \in \orbit} f(x^{-1} \gamma x)
$$
is absolutely integrable over $\AmodZF{G}{G}{Z}$, and its integral is
$$			
	J_\orbit(f) = m(\AmodZF{G}{G}{Z}) \int_{\AmodA{G}{G_\gamma}} f(x^{-1} \gamma x) \d x.
$$
\end{Prop}
\begin{pf} It suffices to prove the integrability assertion, since the ensuing
formula for $J_\orbit(f)$ follows just as in the case of compact quotient.

First note that for $x$ lying in any fixed set $\Omega_1$, which is compact
modulo $\Apts{Z}$, the sum defining $K_\orbit(x, x)$ is actually finite
and hence $K_\orbit$ defines a nice, smooth function on $\AmodZF{G}{G}{Z}$.
Indeed, $f(x^{-1} \gamma x) \ne 0$ implies
$x^{-1} \gamma x \in \Omega_f = $ support of $f$. Thus $\gamma$ itself
must lie in a compact set modulo $\Apts{Z}$, namely
$\Omega_1\Omega_f\Omega_1^{-1}$, and such a set can contain only finitely
many $\gamma$ in $\FmodF{G}{Z}$.

On the other hand, if $\gamma$ is {\em elliptic\/}, then $f(x^{-1} \gamma x) \ne 0$
implies (by Lemma \ref{L:compacity}) that $\alpha H(g) < d_{\Omega_f}$.
{I.}{e.},\ the support of
$$			
      K_\ellip(x, x) = \sum_{\gamma \text{ elliptic}} f(x^{-1} \gamma x)
$$
lies in the ``doubly truncated" Siegel set
$$			
      	\Siegel{c}^d = \Apts{Z} C_1 \{h_t\} C_2 K \word{with} c < 2t < d_{\Omega_f}.
$$
Equivalently, the support of $K_\ellip$ is compact modulo $\Apts{Z}$,
and the integrability is clear.
\end{pf}

\begin{Rem} The argument above actually proves that the {\em full\/} elliptic
contribution is absolutely integrable, \ie, that
$$			
      	\int \abs {K_{\ellip}(x, x)} \d x < \infty.
$$
Thus, by the Lebesgue dominated convergence theorem, it follows that
$$			
      	\sum_{\orbit \text{ elliptic}} \int \abs {K_\orbit(x, x)} \d x < \infty,
$$
and hence that the sum
$$			
      	\sum_{\orbit \text{ elliptic}} J_\orbit(f)
$$
converges absolutely.
\end{Rem}

It remains now to prove Lemma \ref{L:compacity}. It might be tempting to do
so for $\GLtwo$ by simply trying to multiply out some $2 \times 2$
matrices (as we did in \cite[p~201]{sG75}). However, this more quickly
than not leads to an {\em incomplete\/} proof (as happened in fact in \cite{sG75}),
whereas the approach of \cite[p~362]{jA74} works simply and generally, as
we shall show below.

\begin{pf*}{Proof of Lemma \ref{L:compacity}} Let
$$			
      	\rho : G \to G\ell(V) 
$$
denote the adjoint \rep\ of $G$ on the space of {\em trace zero\/} $2 \times 2 $
matrices over $\Nfield$; its highest weight vector is $\Lambda = \alpha$,
the unique simple root of $G$ with respect to the torus $M$. Fix a basis
$\{\e_0, \e_1, \e_2\}$ of $V$ such that
\begin{alignat*}{2}
      	\rho(a) \e_j &= a^{\Lambda - j\alpha}  \e_j
      			 &&\word{for}  a \in M 	\quad (j = 0, 1, 2) 		\\
      	\rho(n) \e_0 &= \e_0  &&\word{for} n \in N \word{and} 	\\
      	\rho(w) \e_0 &= \e_2 	&&
\end{alignat*}
For $\xi$ in $\Apts{V}$, define
$$			
      	\norm{\xi} = \prod_\place \norm{\xi_\place}_\place,
$$
where, for finite $\place$,
$\norm{\xi_\place}_\place = \max_j \{ \abs{\xi_\place^j}_\place\} $
(where $(\xi_\place^j)$ denote the coordinates of $\xi_\place$ with respect
to the basis $\{\e_j\}$), and for infinite $\place$, $\norm{\xi_\place}_\place$
is defined by the Hilbert space structure that makes $\{\e_j\}$ into an \ONB.
(To make sense out of the infinite product, we restrict $\xi$ to be
``primitive", \ie, $\norm{\xi_\place}_\place = 1$ for all but finitely
many $\place$.)

To prove the Lemma, suppose that
\begin{equation} \tag{${*}$}
      	g^{-1} \gamma n g \in \Omega
\end{equation}
 with $n \in N(\Adele)$ and $\Omega$ compact modulo $Z(\Adele)$.
 Write $g = n_1 a k$, so that $({*})$ implies
$$			
      	a^{-1} n_1^{-1} \gamma n n_1 a \in K \Omega K.
$$
Because the map $g \mapsto \rho(g) \e_0$ is continuous, and $\rho$ is
trivial on $\Apts{Z}$, it follows that
$$			
      	\norm{\rho(a^{-1} n_1^{-1} \gamma n n_1 a )\cdot \e_0} \le \e^{2 d_0}
$$
for some $d_0 > 0$. So suppose that $\gamma$ does not belong to
$B(\Nfield)$. Then by Bruhat's decomposition,
$$			
      	\gamma = b_0 w n_0
$$
for some  $b_0$ in $B(\Nfield)$ and $n_0$ in $N(\Nfield)$. But for any
$b \in \Apts{B}$ and $g \in \Apts{G}$, our assumptions on $\rho$ and
$\e_0$ imply
$$			
      	\norm{\rho(g b)\cdot \e_0} = \e^{\alpha H(b)}\, \norm{\rho(g)\cdot \e_0}.
$$
Thus (with $\gamma = b_0 w n_0$) we compute that
\begin{align*}
	\norm{\rho((a^{-1} n_1^{-1} b_0 w)  (n_0 n n_1 a ))\cdot \e_0}
		&= \e^{\alpha H(a)}\, \norm{\rho(a^{-1} n_1^{-1} b_0 )\cdot \e_2}	\\
		&= \e^{ 2\alpha H(a)}\, \norm{\rho( n^*)\cdot \e_2}			\\
		&\ge \e^{ 2\alpha H(a)}\, \norm{\e_2} = \e^{ 2\alpha H(g)}
\end{align*}
(since $H(g) = H(a)$, $\rho(a )\, \e_2 = a^{-\alpha} \,\e_2 $, and
$\norm{\rho( n) \,\e_2} \ge \norm{\e_2}$ for any $n$ in $\Apts{N}$).
Thus the Lemma is proved for $d_\Omega  > d_0$.
\end{pf*}

\bigskip
\subsection{The Definition of \protect$\bold{K_\orbit^T(x, y)}\protect$} \label{S:sec3}

For $\orbit$ the (unique non-trivial) unipotent class in $\FmodF{G}{Z}$, we saw
in Lecture I that $ K_\orbit(x, y)$ is {\em not\/} integrable over the diagonal. A
somewhat less trivial computation  (\cf Lecture III)
shows that the hyperbolic $K_\orbit(x, y)$
are also not integrable. So let us finally introduce the required modifications of
these kernels.

For $P = B$, let $\trunc B$ denote the characteristic function of the
{\em positive Weyl chamber\/}
$$			
      	\WeylCh B = \set{(r_1, r_2) \in \LAa B = \LAa M}{r_1 - r_2 > 0},
$$
\ie, the cone in $\R^2$ where $\alpha(r_1, r_2) = r_1 - r_2$
is positive. Then define, for any $\orbit \in \Orbit$ and $T = (T_1, T_2)$
in $\WeylCh B$,
\begin{align} \label{E:modifiedKer}
 	\GeoKer(x, f) &= K_\orbit(x, x) - \sum_{\delta \in \FmodF{G}{B}}
 				   K_{B, \orbit}(\delta x, \delta x)\ \trunc B(H(\delta x) - T)	\\
 	\intertext{where}
 	K_{B, \orbit}(x, y) &= \sum_{\gamma \in \orbit \cap \FmodF{M}{Z}}
 					\int_{\Apts{N}} f(x^{-1} \gamma n y) \d n. \notag
 \end{align}
This is Arthur's {\em modified kernel}\ \ $\GeoKer(x, f)$.

Note that if $\orbit$ is elliptic, then $\orbit \cap M(\Nfield) = \emptyset$.
In this case $K_{B, \orbit} \equiv 0$, and
$$			
      	\GeoKer( x, f) = K_\orbit(x, x).
$$
In general, the term subtracted in \ref{E:modifiedKer} from $K_\orbit(x, x)$ is
the ``correction term at $\infty$'' which ends up making $\GeoKer(x, f)$
integrable. Before proving this is so, let us note the inductive nature of the
definition of $\GeoKer$ (thereby making more plausible the appearance of
$K_{B, \orbit}$). Set
\begin{align} \label{E:KB}
	K_B(x, y)	&= \sum_{\orbit} K_{B, \orbit}(x, y)	\\
 			&\equiv \sum_{\gamma \in  \FmodF{M}{Z}}
 					\int_{\Apts{N}} f(x^{-1} \gamma n y) \d n.  \notag
 \end{align}
 Then $K_B(x, y)$ is just the kernel of $R(f)$ acting in
 $L^2(\EisenSp{B}{N})$ in place of $L^2(\AmodZF{G}{G}{Z}$.
 Indeed for $\varphi$ in $L^2(\EisenSp{B}{N})$,
\begin{align} \label{E:Rf}
      	R(f)\varphi(x) &= \int_{\AmodA{G}{Z}} f(x^{-1} y) \varphi(y) \d y	 \notag	\\
      				&= \int_{\EisenSp{B}{N}} \biggl( \sum_{\FmodF{M}{Z}}
      				     \int_{\Apts{N}} f(x^{-1} \gamma n y) \d n \biggr)
      				     \varphi(y) \d y,
 \end{align}
 as claimed. In general, for an arbitrary $G$, the definition of $\GeoKer(x, f)$
 takes into account {\em all\/} the proper parabolics $P$ of $G$, subtracting
 off from $K_\orbit(x, x)$ not just one term (for the minimal parabolic $B$),
 but the sum
 \begin{align} \label{E:Psum}
 	&\sum_{P \subset G} (-1)^{\dim A_P / Z} \sum_{\delta \in \FmodF{G}{P}}
 				   K_{P, \orbit}(\delta x, \delta x)\ \trunc P (H(\delta x) - T),	\\
 	\intertext{where}
 	K_P(x, y) &= \sum_{\orbit \in \Orbit} K_{P, \orbit}(x, y)   \notag
 \end{align}
is the kernel of $R_P(f)$ in $L^2(\EisenSp{M_P}{N_P})$.  (Here $A_P$ is the
maximal split torus in $M_P$, the Levi component of $P = M_P N_P$, and $R_P$
is the right regular \rep\ of $\Apts{G}$ in $L^2(\EisenSp{M_P}{N_P})$.)

\begin{Rem}
{\rm
Observe that
$$			
      	\GeoKer( x, f) = K_\orbit(x, x).
$$
for $x$ in a compact (modulo $\Apts{Z}$) set $\Omega$ (how large depends
on T). Indeed, for $x$ in an appropriate such set $\Omega$,
$\trunc B(H(\delta x) - T)$ is identically
zero\footnote{{\em Proof:\ } Fix $c$				
and $T$ such that $\alpha(T) > c$. Fact \ref{F:FundDom1} of Section \ref{S:sec2}
implies that in $\Siegel c$, $\trunc{B}(H(\delta x) - T) \ne 0$ only for
finitely many classes $\delta_i$ in $\FmodF{G}{B}$. Then for each such class,
$\trunc{B}(H(\delta_i x) - T) = 0$ on the compact set $\Omega_i^T$ in
$\Siegel c$ where $c < H(\delta_i x) < \alpha(T)$. Thus
$\trunc{B}(H(\delta x) - T)$ is identically zero on
$\Omega_T = \cap_{i = 1}^n \Omega_i^T$.}.			
Thus $x$ is indeed a modification of $K_\orbit(x, x)$ only ``near infinity''
(and in higher rank, more than one parabolic is needed to effect these
modifications near infinity \dots).
}
\end{Rem}

\bigskip
\subsection{The Integrability of $\GeoKer( x, f)$} \label{S:sec4}

 Henceforth, let us assume that $\orbit$ is unipotent or hyperbolic. To prove that
 $\GeoKer( x, f)$ is integrable over $\AmodZF{G}{G}{Z}$, it suffices to
 show that it is integrable over a Siegel domain $\Siegel c$. So pick any
 $T \in \WeylCh B$, and write
\begin{align*}
      	\Siegel c &= \SiegelPiece{T} \cup \SiegelPieceT{T},	\\
      	\intertext{where}
      	\SiegelPiece{T} &= \set{g \in \Siegel c}{\alpha H(g) > \alpha(T) > 0}, 	
\end{align*}
and $ \SiegelPieceT{T}$ is its complement in $\Siegel c$. Clearly
$\GeoKer(x, f)$ is integrable over $\SiegelPieceT{T}$,  since
$$			
	c < \alpha H(g) \le \alpha(T).
$$
Thus it suffices to prove that  $\GeoKer( x, f)$
is integrable over $\SiegelPiece{T}$. (\NB In general, the partition of a
Siegel domain requires more non-trivial ``reduction theory'' than is evident in
the case of $\GLtwo$; this is what leads to the ``geometric'' combinatorial
lemmas in \S 6 of \cite{jA78}.)
\begin{Lem} \label{L:BddAlpha}
There exists a constant $c_0$ such that $\alpha H(wn) \le c_0$
for all $n \in \Apts{N}${\rm ;} moreover, if $\gamma$ in $G(\Nfield)$
satisfies the inequality
$$			
	\alpha H(\gamma x) > c_0
$$
for some $x$ in $\Apts{G}$ with $\alpha H(x) > c_0$, it must follow that
$\gamma \in B(\Nfield)$.
\end{Lem}
As the proof is similar to that of the crucial Lemma \ref{L:compacity} in
\S \ref{S:sec2} of this Lecture, we leave the details to the reader (\cf
\cite[pp~334--335]{jA74}).

Now {\em fix} $T$ in $\WeylCh B$ such that $\alpha (T) > c_0$. Then for $x$
in $\SiegelPiece{T}$, Lemma \ref{L:BddAlpha} implies that
$$			
	\trunc{B}(H(\delta x) - T) \ne 0
$$
(if and) only if $\delta \in B(\Nfield)$. Thus (for such $x$),
 \begin{align*}
 	\sum_{\gamma \in \orbit} f(x^{-1} \gamma x) =
 	\sum_{\delta \in \FmodF{G}{B}} &\sum_{\gamma \in \orbit}
 	f(x^{-1} \delta^{-1} \gamma \delta x) \trunc{B}(H(\delta x) - T),   \\
 	\intertext{and hence}
 	\GeoKer(x, f) = \sum_{\delta \in \FmodF{G}{B}}
 	\biggl\{ \sum_{\gamma \in \orbit}\, &f(x^{-1} \delta^{-1} \gamma \delta x) 	\\
 	- \int_{\Apts{N}} &\sum_{\mu \in M(\Nfield) \cap \orbit}
 	f(x^{-1} \delta^{-1} \mu n \delta x) \d n \biggr\} \trunc{B}(H(\delta x) - T)	
 \end{align*}

Next fix $T$ such that $\alpha(T) > \max\{c_0,  d_{\Omega_f}\}$, with
$ d_{\Omega_f}$ as in Lemma~\ref{L:compacity}. Then by
Lemma~\ref{L:compacity}, we may replace the sum here over
$\gamma \in \orbit$ by a sum over
$$						
	B(\Nfield) \cap \orbit = N(\Nfield) \cdot M(\Nfield) \cap \orbit
$$
(This last equality holds since the class of $\gamma$ in $B(\Nfield) $ is
determined by its semi-simple part, \ie its $M(\Nfield)$-part.) Thus we get

\begin{align} \label{E:Poisson}
	\GeoKer( x, f) &=
	\sum_{\delta \in \FmodF{G}{B}} \sum_{\mu \in M(\Nfield) \cap \orbit}
 	\biggl\{ \sum_{\nu \in N(\Nfield)} f(x^{-1} \delta^{-1} \mu\nu \delta x) \notag\\	
 	& \qquad \qquad \qquad  \qquad \qquad
 	- \int_{\Apts{N}} f(x^{-1} \delta^{-1} \mu n \delta x) \d n \biggr\}
 	\trunc{B}(H(\delta x) - T)										  	\\
 	&= \sum_\delta \sum_\mu \sum_{\nu \in N(\Nfield) - \{1\}}
 	\widehat{F}_{\delta x, \mu}  (\nu) \, \trunc{B}(H(\delta x) - T)		    \notag\\
 	 \intertext{ where }
	& \qquad \qquad {F}_{ x, \mu} (n) = f(x^{-1} \mu n x)				    \notag \\
 	 \intertext{defines a Schwartz-Bruhat function on $\Apts{N} \approx \Adele$
			   with Fourier transform }
 	 & \qquad \qquad \widehat{F}_{x, \mu} (\nu)  =
 	 		\int_{\Apts{N}} f(x^{-1} \mu n x) \psi_N(n \nu) \d n	\notag 	
\end{align}
(Here $\psi_N(\bb 1 x 1) = \psi(x)$ for some fixed non-trivial character on
$\Nfield \bs \Adele$, and we have applied the Poisson summation formula to
${F}_{\delta x, \mu}$.) Finally we are ready to prove

\begin{Prop}  \label{P:integrable} For any
$\orbit \in \Orbit$, $\GeoKer(x, f)$
is absolutely integrable over 
$$\AmodZF{G}{G}{Z}.$$
\end{Prop}

\begin{pf}
We may assume $\orbit$ is not elliptic.
By the above arguments, it suffices to check the integrability over
$\SiegelPiece T$, with $\alpha(T) > \max\{c_0,  d_{\Omega_f}\}$ as above.
As in that domain formula \eqref{E:Poisson} holds, we conclude
\begin{align} \label{E:IntAbs}
	 \int \abs{\GeoKer( x, f)} \d & x  \le							 \\
	 \int_{\AmodZF GBZ} &\sum_{\mu \in M(\Nfield) \cap \orbit}
	 \sum_{\nu \in N(\Nfield) - \{1\}}
 	\abs{\widehat{F}_{x, \mu} (\nu)} \,\trunc{B}(H(x) - T) \d x.  \notag
\end{align}
Note that for any  $\mu = {\bb \alpha 0 1}$ in $M(\Nfield) \cap \orbit$,
$\nu = {\bb 1 w 1}$ in $N(\Nfield) - \{1\}$, and $x = z n h_t m k$,
\begin{align*}
	\widehat{F}_{x, \mu} (\nu) &=
	\int_\Adele f\left( x^{-1} {\bb \alpha 0 1} {\bb 1 y 1} x \right) \psi(y w) \d y 	\\
	&= \e^{2t} \widehat{F}_{k^*, \mu} (\e^{2t} \nu)
\end{align*}
with $k^* = (h_t^{-1} n h_t) m k$. (Here $\e^{2t} \nu = \bb {1}{\e^{2t} w}{1}$,
and we confuse $N(\Nfield)$  with $\Nfield$.) Note also that the variable
$x$ in \eqref{E:IntAbs} runs through a Siegel set. Thus $k^*$ runs through a
compact set, and we have, for any $N > 0$
$$			
	\sum_{\nu \ne 1} \widehat{F}_{k^*, \mu} (\e^{2t} \nu) \le c_N (\e^{-2t})^N
$$
Applying Iwasawa's decomposition \eqref{E:IwasawaInt} to \eqref{E:IntAbs},
we then conclude that
$$			
	\int \abs{\GeoKer( x, f)} \d x \le
	C \int_{\AmodF{N}} \int_{\Nfield^\times \bs \Adele^1} \int_K
	\biggl(\int_{\alpha(T)}^\infty \e^{-2tN} \d t \biggr) \d k \d^\times a \d n < \infty
$$
\nopagebreak as required.
\end{pf}

Summing up, we have
\begin{Thm} \label{T:AbsInt}
	For $\alpha(T)$ sufficiently large,
$$						
\sum_{\orbit \in \Orbit} \int_{\AmodZF{G}{G}{Z}}
						\abs{\GeoKer( x, f)} \d  x  < \infty. \qed
$$
\end{Thm}
Hence the left side $\sum_{\orbit \in \Orbit} J^T_\orbit (f)$ of
(1) of Lecture I 
is defined and absolutely convergent.

Indeed, it suffices to remark that for a given $f \in \Cinfc(\AmodA{G}{Z})$,
$\GeoKer( x, f)$  is identically zero except for finitely many hyperbolic
$\orbit$. This is clear from the defining equation \eqref{E:modifiedKer},
and the fact that $x^{-1} \delta^{-1} {\bb \alpha 0 1}\delta x$ and
$x^{-1} \delta^{-1} {\bb \alpha 0 1} n \delta x$ can belong to the fixed
compact mod $\Apts{Z}$ set $\Omega_f$ for only finitely many $\alpha$;
to see this, observe that conjugating by $\delta x$ preserves
 eigenvalues.  Thus the sum over
$\orbit$ in $\Orbit$ in Theorem \ref{T:AbsInt} is infinite only over the
elliptic classes, where we already know the relevant sum 
$\sum \int \abs{K_\orbit (x, f)} \d x $ is finite  (\cf the Remark
after the proof of Proposition \ref{P:elliptic}).

\bigskip\noindent
{\bf Concluding Remarks.} (a) The terms $J^T_\orbit (f)$ are indeed
generalizations of the orbital integrals appearing on the geometric side of
Selberg's compact quotient trace formula. Indeed, for $\orbit$ elliptic
(which is automatic for $G$ anisotropic), $J^T_\orbit (f)$ reduces to a
multiple of the usual orbital integral
$\int_{\AmodA{G}{G_\gamma}} f(x^{-1} \gamma x) \d x$.

\medskip\noindent
We have seen---for $\GLtwo$---that the kernel $K_\orbit (x, x)$
is already integrable on $\Siegel c$ if $\orbit$ is an {\em elliptic class}
(as opposed to hyperbolic or unipotent). For general $G$, the
proper distinction is between those $\orbit$ which never intersect
{any proper} parabolic subgroup subgroup $P$ of $G$ and those that do.
Indeed Arthur does not talk about ``elliptic'' $\orbit$ at all, but rather defines
a modified kernel
$$			
	\GeoKer( x, f)
$$
for any $\orbit$, and proves that $\GeoKer( x, f)$ is integrable. But of
course---as we have already observed---$\GeoKer( x, f) = K_\orbit(x, x)$
precisely when $\orbit$ is elliptic in this sense (that it intersects no
proper $P$ \dots).

\end{document}